% Plain TeX

%\magnification=1050
%\hsize=6.5truein
%\vsize=8.25truein
%\hoffset=1.0truein
%\voffset=0.5truein
%\baselineskip=20truept
%\baselineskip=24pt % double linespace

\input eplain

\newcount\fignumber
\def\figdef#1{\global\advance\fignumber by 1 \definexref{#1}{\number\fignumber}{figure}\ref{#1}}
\def\figdefn#1{\global\advance\fignumber by 1 \definexref{#1}{\number\fignumber}{figure}}
\let\figref=\ref
\let\figrefn=\refn
\let\figrefs=\refs

\newcount\tabnumber
\def\tabdef#1{\global\advance\tabnumber by 1 \definexref{#1}{\number\tabnumber}{table}\ref{#1}}
\def\tabdefn#1{\global\advance\tabnumber by 1 \definexref{#1}{\number\tabnumber}{table}}

% Thanks to http://insti.physics.sunysb.edu/~siegel/tex.shtml
%
% for ordinary tex
\ifx\pdfoutput\undefined
\input epsf

\def\figscale#1#2{\epsfxsize=#2\epsfbox{#1.eps}}
%
% for pdftex
\else

\def\figscale#1#2{\pdfximage width#2 {#1.pdf}\pdfrefximage\pdflastximage}
\fi

%% from /usr/share/texlive/texmf-dist/tex/plain/base/gkpmac.tex
%% \def\beginsection#1 #2 \par{ % should not be first in the chapter
%%   \backup=\lastskip	% but should come on first or second page of chapter
%%   \mark{#1\enspace #2}	% because the mark gives running head on right page
%%   \nobreak\vskip-\backup\penalty-200
%%   {\subtitle\baselineskip=34pt
%%     \noindent\hbox to2\parindent{#1\hfil}\uppercase{\kern-.05em#2}\par}
%%   \nobreak\vskip5pt\noindent\hbox to2\parindent{}}

\newcount\scount \scount=0

%% \def\section#1{
%%     \global\advance\scount by1
%%     \definexref{#1}{\the\scount}{section}
%%     \vskip.25truein\noindent\the\scount\quad{\bf #1}\hfill\vskip.25truein
%%     }

%% from /usr/share/texlive/texmf-dist/tex/plain/base/plain.tex
%% \outer\def\beginsection#1\par{\vskip\z@ plus.3\vsize\penalty-250
%%   \vskip\z@ plus-.3\vsize\bigskip\vskip\parskip
%%   \message{#1}\leftline{\bf#1}\nobreak\smallskip\noindent}

\makeatletter
\def\section#1\par{
  \vskip\z@ plus.3\vsize\penalty-250
  \vskip\z@ plus-.3\vsize\bigskip\vskip\parskip
  \global\advance\scount by1
  \writenumberedtocentry{section}#1{}
  \definexref#1{\the\scount}{section}
  \message{#1}
  \noindent\the\scount.\quad{\bf #1}\nobreak\smallskip\noindent}
\makeatother

\centerline{\bf{Easy Accurate Reading and Writing of Floating-Point Numbers}}
\bigskip
\centerline{Aubrey Jaffer\numberedfootnote{Digilant, 2 Oliver Street Suite 901, Boston, MA 02109.  Email: {\tt agj@alum.mit.edu}}}
\bigskip
\centerline{August 2018}
\bigskip

\beginsection{Abstract}

{\narrower
  Presented here are algorithms for converting between (decimal)
  scientific-notation and (binary) IEEE-754 double-precision
  floating-point numbers.  By employing a rounding integer quotient
  operation these algorithms are much simpler than those previously
  published.  The values are stable under repeated conversions between
  the formats.  Unlike Java-1.8, the scientific representations
  generated use only the minimum number of mantissa digits needed to
  convert back to the original binary values.

  %% Implemented in Java these algorithms execute as fast or faster than
  %% Java's native conversions over nearly all of the IEEE-754
  %% double-precision range.
\par}

%% \beginsection{Table of Contents}

%% \item{1} Introduction
%% \itemitem{1.1} Relation to Previous Work
%% \itemitem{1.2} Symbols
%% \item{2} BigIntegers
%% \item{3} Reading
%% \item{4} Writing
%% \item{5} Performance
%% \item{6} Conclusion
%% \item{} References

\beginsection{Introduction}

  Articles from Steele and White\cite{Steele:1990:PFN:93548.93559},
  Clinger\cite{Clinger:1990:RFP:93548.93557}, and Burger and
  Dybvig\cite{Burger:1996:PFN:249069.231397} establish that binary
  floating-point numbers can be converted into and out of decimal
  representations without losing accuracy while using a minimum number
  of (decimal) significant digits.  Using the minimum number of digits
  is a property which Java-1.8 does not achieve ($10^{23}$ prints as
  {\tt 9.999999999999999E22}; $8\times10^{-323}$ prints as {\tt
  7.9E-323}); the {\tt doubleToString} procedure presented here
  produces only minimal precision mantissas.

  The lossless algorithms from these papers all require high-precision
  integer calculations, although not for every conversion.

  In {\it How to Read Floating-Point Numbers
  Accurately}\cite{Clinger:1990:RFP:93548.93557} Clinger astutely
  observes that successive rounding operations do not have the same
  effect as a single rounding operation.  This is the crux of the
  difficulty with both reading and writing floating-point numbers.
  But instead of constructing his algorithm to do a single rounding
  operation, Clinger and the other authors follow
  Matula\cite{Matula:1968:IC:362851.362887,1671610} in doing
  successive roundings while tracking error bands.

  The algorithms from {\it How to Print Floating-point Numbers
  Accurately}\cite{Steele:1990:PFN:93548.93559} and {\it Printing
  floating-point numbers quickly and
  accurately}\cite{Burger:1996:PFN:249069.231397} are iterative and
  complicated.  The read and write algorithms presented here do at
  most 2 and 4 BigInteger divisions,
  respectively\numberedfootnote{Writing exact powers of two takes up
  to 6 BigInteger divisions, but there are only 2100 of them in the
  IEEE-754 double-precision range; they could be precomputed.}
  %%   This simplicity is responsible for the speed of the algorithms.

  Over the range of IEEE-754\cite{IEEE:1985:AIS} double-precision
  numbers, the largest intermediate BigInteger used by these
  power-of-5 algorithms is 242 decimal digits (803 bits).  Steele and
  White\cite{Steele:1990:PFN:93548.93559} report that the largest
  integer used by their algorithm is 1050 bits.  These are not large
  for BigIntegers, being orders of magnitude smaller than the smallest
  precisions which get speed benefits from FFT multiplication.

  Both Steel and White\cite{Steele:1990:PFN:93548.93559} and
  Clinger\cite{Clinger:1990:RFP:93548.93557} claim that the input and
  output problems are fundamentally different from each other because
  the floating-point format has a fixed precision while the decimal
  representation does not.  Yet, in the algorithms presented here,
  BigInteger rounding divisions accomplish %% fast
  accurate conversions in both directions.

  While the read algorithm tries the division yielding the longer
  precision quotient first, and retries only if it doesn't fit into
  the mantissa, the write algorithm tries the shorter precision
  division first and retries only when the shorter precision fails to
  read back correctly.

\beginsection{BigIntegers}

  Both reading and writing of floating-point numbers can involve
  division of numbers larger than can be stored in the floating-point
  registers, causing rounding at unintended steps during the
  conversion.

  BigIntegers (arbitrary precision integers) can perform division of
  large integers without rounding.  What is needed is a BigInteger
  division-with-rounding operator, called {\tt roundQuotient} here.
  For positive operands, it can be implementated in Java as follows:

\smallskip\verbatim|
public static BigInteger roundQuotient(BigInteger num, BigInteger den) {
    BigInteger quorem[] = num.divideAndRemainder(den);
    int cmpflg = quorem[1].shiftLeft(1).compareTo(den);
    if (quorem[0].and(BigInteger.ONE).equals(BigInteger.ZERO) ?
        1==cmpflg : -1<cmpflg)
        return quorem[0].add(BigInteger.ONE);
    else return quorem[0];
}
|endverbatim
\medskip

  If the remainder is more than half of the denominator, then it
  rounds up; if it is less, then it rounds down; if it is equal, then
  it rounds to even.  These are the same rounding rules as the IEEE
  Standard for Binary Floating-Point Arithmetic\cite{IEEE:1985:AIS}.

  For the algorithms described here the value returned by
  roundQuotient always fits within a Java long.  Having roundQuotient
  return a Java {\tt long} integer turns out to execute more quickly
  than when a BigInteger is returned.

\smallskip\verbatim|
public static long roundQuotient(BigInteger num, BigInteger den) {
    BigInteger quorem[] = num.divideAndRemainder(den);
    long quo = quorem[0].longValue();
    int cmpflg = quorem[1].shiftLeft(1).compareTo(den);
    if ((quo & 1L) == 0L ? 1==cmpflg : -1<cmpflg) return quo + 1L;
    else return quo;
}
|endverbatim
\medskip

  %% For its floating-point conversions Java uses a small special-purpose
  %% big-integer implementation named {\tt FDBigInt}.  The
  %% implementations of the new algorithms presented here use the Java
  %% {\tt BigInteger} package, which comes with the Java distribution.

  In the scaled twos-complement encoding of the mantissa, the
  representation of 5 and 10 are the same; the exponent differs by
  one.  The same is true of any non-negative integer power of 5 and
  10.

  In the algorithms below, {\tt bipows5} is an array of 326 BigInteger
  successive integer powers of 5.  Constant {\tt dblMantDig} is the
  number of bits in the mantissa of the normalized floating-point
  format (53 for IEEE-754 double-precision numbers).  Constant {\tt
  llog2} is the base 10 logarithm of 2.

\beginsection{Reading}

  The {\tt MantExpToDouble} algorithm computes the closest (binary)
  floating-point number to a given number in scientific notation by
  finding the power-of-2 scaling factor which, when combined with the
  power-of-10 scaling specified in the input, yields a
  rounded-quotient integer which just fits in the binary mantissa
  (having {\tt dblMantDig} bits).

  The first argument, {\tt lmant}, is the integer representing the
  string of mantissa digits with the decimal point removed.  The
  second argument, {\tt point}, is the (decimal) exponent less the
  number of digits of mantissa to the left of the decimal point.  If
  there was no decimal point it is treated as though it appears to the
  right of the least significant digit.  Thus {\tt point} will be zero
  when the floating-point value equals the integer {\tt lmant}.

  When {\tt point} is non-negative, the mantissa is multiplied by
  $5^{\rm\bf point}$ and held in variable {\tt num}.  If {\tt num}
  fits within the binary mantissa, {\tt num.doubleValue()} converts to
  the correct double-precision mantissa value and {\tt Math.scalb}
  scales by {\tt point} bits.  Otherwise {\tt MantExpToDouble} calls
  {\tt roundQuotient} to divide and round to {\tt dblMantDig} bits.
  Because the divisor is a power-of-2, the number of bits in the
  quotient is one more than the difference of the number of bits of
  dividend and divisor.

  With a negative {\tt point}, the mantissa will be multiplied by a
  power of 2, then divided by {\tt scl} $=5^{\rm\bf - point}$.  To scale
  by $2^{\rm\bf - point}$, {\tt point} is added to {\tt bex} to form the
  binary exponent for the returned floating-point number.

  The integer quotient of a $n$-bit positive integer and a smaller
  $m$-bit positive integer ($0<m<n$) will always be between $n-m$ and
  $1+n-m$ bits in length.  Because rounding can cause a carry to
  propagate through the quotient, the longest integer returned by the
  {\tt roundQuotient} of a $n$-bit positive integer and a smaller
  $m$-bit positive integer is $2+n-m$ bits in length, for example {\tt
  roundQuotient}$(7,2)\to4$.  If this happens for some power-of-five
  divisor (which is close to a power of 2) then it must happen when
  the dividend is the largest possible $n$-bit integer, $2^n-1$.

  Over the double-precision floating-point range (including
  denormalized numbers) there are only 2100 distinct positive numbers
  with mantissa values which are all (binary) ones ($2^n-1$); testing
  all of them finds that in doing double-precision floating-point
  conversions, there is no integer power-of-5 close enough to an
  integer power-of-2 which, as divisor, causes the quotient to be
  $2+n-m$ bits in length.

  Thus the longest a rounded-quotient of a $n$ bit integer and a $m$
  bit power-of-5 can be is $1+n-m$ bits; the shortest is $n-m$ bits.
  This means that no more than 2 rounded-quotients need be computed in
  order to yield a mantissa which is {\tt mantlen} bits in length.

\verbatim|

public static double MantExpToDouble(long lmant, int point) {
    BigInteger mant = BigInteger.valueOf(lmant);
    if (point >= 0) {
        BigInteger num = mant.multiply(bipows5[point]);
        int bex = num.bitLength() - dblMantDig;
        if (bex <= 0) return Math.scalb(num.doubleValue(), point);
        long quo = roundQuotient(num, BigInteger.ONE.shiftLeft(bex));
        return Math.scalb((double)quo, bex + point);
    }
    BigInteger scl = bipows5[-point];
    int mantlen = dblMantDig;
    int bex = mant.bitLength() - scl.bitLength() - mantlen;
    int tmp = bex + point + 1021 + mantlen;
    if (tmp < 0) {bex -= tmp + 1; mantlen += tmp;}
    BigInteger num = mant.shiftLeft(-bex);
    long quo = roundQuotient(num, scl);
    if (64 - Long.numberOfLeadingZeros(quo) > mantlen)
        {bex++; quo = roundQuotient(num, scl.shiftLeft(1));}
    return Math.scalb((double)quo, bex + point);
}
|endverbatim
%% \vfill\eject

  The lines involving {\tt tmp} reduce {\tt mantlen} for denormalized
  floating-point representation when the number is too small for the
  floating-point exponent.  {\tt bex} and {\tt mantlen} are offset by
  different amounts because the normalized mantissa has an implied
  most significant 1 digit, while it is explicit in denormalized
  mantissas.

  When {\tt point}$<0$, if the number returned by the call to {\tt
  roundQuotient} is more than {\tt mantlen} bits long, then call
  {\tt roundQuotient} with double the denominator {\tt scl}.  In
  either case, the final step is to convert to floating-point and
  scale it using {\tt Math.scalb}.

  Because the quotient which gets used is rounded by a single
  operation to the correct number of bits, it is the closest to the
  decimal value possible in binary floating-point representation.

  Separating powers of 2 from powers of 5 in the {\tt MantExpToDouble}
  algorithm enables a 29\% reduction in the length of intermediate
  BigIntegers.

  In {\it Fast Path Decimal to Floating-Point
  Conversion}\cite{Regan2011fastpath} Regan describes using
  floating-point multiplication and division when the mantissa and
  power-of-ten scale fit within floating-point mantissas.  This second
  version of {\tt MantExpToDouble} uses floating-point for small
  magnitude expoenents, separating the powers of 5 and powers of 2 and
  postponing the binary scaling until after the multiplication or
  division by power of 5, which extends the range for which Regan's
  method applies.

\verbatim|
public static double MantExpToDouble(long lmant, int point) {
    long quo; int bex;
    if (point >= 0) {
        if (point < dpows5.length &&
            64 - Long.numberOfLeadingZeros(lmant) <= dblMantDig)
            return Math.scalb(((double)lmant) * dpows5[point], point);
        BigInteger mant = BigInteger.valueOf(lmant);
        BigInteger num = mant.multiply(bipows5[point]);
        bex = num.bitLength() - dblMantDig;
        quo = roundQuotient(num, BigInteger.ONE.shiftLeft(bex));
        return Math.scalb((double)quo, bex + point);
    }
    if (-point < dpows5.length &&
        64 - Long.numberOfLeadingZeros(lmant) <= dblMantDig)
        return Math.scalb(((double)lmant) / dpows5[-point], point);
    BigInteger mant = BigInteger.valueOf(lmant);
    BigInteger scl = bipows5[-point];
    int mantlen = dblMantDig;
    bex = mant.bitLength() - scl.bitLength() - mantlen;
    int tmp = bex + point + 1021 + mantlen;
    if (tmp < 0) {bex -= tmp + 1; mantlen += tmp;}
    BigInteger num = mant.shiftLeft(-bex);
    quo = roundQuotient(num, scl);
    if (64 - Long.numberOfLeadingZeros(quo) > mantlen)
        {bex++; quo = roundQuotient(num, scl.shiftLeft(1));}
    return Math.scalb((double)quo, bex + point);
}
|endverbatim
\vfill
\beginsection{Writing}

  The goal for writing a floating-point number is to output the
  shortest decimal mantissa which reads back as the original
  floating-point input.  But there are subtleties to this simple
  sounding idea.

  Consider reading back a power-of-two; the number of bits read back
  can depend on which way the decimal output was rounded.  If it was
  rounded up, then the binary mantissa will have the correct number of
  bits.  If it was rounded down, then the binary mantissa will be one
  bit short; if that decimal representation doesn't read back
  correctly, then the binary-to-decimal algorithm recomputes with more
  precision.

  In some cases, the rounded up number reads back correctly, even
  though the rounded down number is more accurate (but with one more
  decimal digit).  This is not the best idea when reading a number
  into a higher precision floating-point format than the format it was
  written from.  The number read may not be the closest to the
  original value.  Burger and
  Dybvig\cite{Burger:1996:PFN:249069.231397} use a different criteria:
  the decimal number written should be the shortest mantissa correctly
  rounded decimal number.  Both cases are treated below.

  The integer quotient of a $n$-digit positive integer and a smaller
  $m$-digit positive integer ($0<m<n$) will always be between $n-m$
  and $1+n-m$ decimal digits in length.  Because rounding can cause a
  carry to propagate through the quotient, the longest integer
  returned by the {\tt roundQuotient} of a $n$-digit positive integer
  and a smaller $m$-digit positive integer is $2+n-m$ digits in
  length, for example {\tt roundQuotient}$(995,10)\to100$.

  A starved precision is tried first; if that does not read back
  correctly, then the written precision is increased by one decimal
  digit; if that does not read back correctly, then the written
  precision is increased by a second decimal digit.  The only cases
  where starved precision correctly rounded numbers are longer than
  necessary are when the trailing digits are zero.  Code at the end of
  the algorithm truncates strings of least-significant 0 digits.

  It turns out that this second extra digit is needed only for powers
  of two and, for IEEE-754 double precision format, only normalized
  powers of two.  Thus the test for this condition simply compares the
  binary mantissa with the largest power of two possible for the
  mantissa, $2^{53}$.  Before trying the second extra digit, the
  quotient plus 1 is checked whether it reads back correctly:

\verbatim|
                if (MantExpToDouble(++lquo, point) != f) {
|endverbatim

  In order to implement the Burger and Dybvig criteria, replace the
  two occurrences of that line by:

\verbatim|
                {
|endverbatim

  In the algorithm, the positive integer mantissa {\tt mant} and
  integer exponent (of 2) {\tt e2} are extracted from floating-point
  input {\tt f}.  Constant {\tt llog2} is the base 10 logarithm of 2.
  The variable {\tt point} is set to the upper-bound of the decimal
  approximation of {\tt e2}, and would be the output decimal exponent
  if the decimal point were to the right of the mantissa least
  significant digit.

  When {\tt e2} is positive, {\tt point} is the upper-bound of the
  number of decimal digits of {\tt mant} in excess of the
  floating-point mantissa's precision.  {\tt mant} is left shifted by
  {\tt e2} bits into {\tt num}.  The {\tt roundQuotient} of {\tt num}
  and $5^{\rm\bf point}$ yields the integer decimal mantissa {\tt lquo}.
  If {\tt mantExpToDouble(lquo,point)} is not equal to the original
  floating-point value {\tt f}, then the {\tt roundQuotient} is
  recomputed with the divisor effectively divided by 10, yielding one
  more digit of precision.

  When {\tt e2} is negative, {\tt den} is set to $2^{\rm\bf -e2}$ and
  {\tt point} is the negation of the lower-bound of the number of
  decimal digits in {\tt den}.  {\tt num} is bound to the product of
  {\tt mant} and $5^{\rm\bf point}$.  The {\tt roundQuotient} of {\tt
  num} and {\tt den} produces the integer {\tt lquo}.  If {\tt
  mantExpToDouble(lquo,point)} is not equal to the original
  floating-point value {\tt f}, then the {\tt roundQuotient} is
  computed again with {\tt num} multiplied by 10, yielding one more
  digit of precision.

  The last part of {\tt doubleToString} constructs the output using
  Java {\tt StringBuilder}.  The mantissa trailing zeros are
  eliminated by scanning the {\tt sman} string in reverse for non-zero
  digits and the decimal point is shifted to the most significant
  digit.

  The Java code for {\tt doubleToString} shown below uses powers of 5
  instead of 10 for speed.  The arguments to {\tt
  BigInteger.leftShift} are adjusted accordingly to be differences of
  {\tt e2} and {\tt point}.

\vfill\eject
\verbatim|
public static String doubleToString(double f) {
    if (f != f) return "NaN";
    if (f+f==f)
        return 1/f<0?"-0.0":(f==0.0)?"0.0":((f > 0) ? "Infinity" : "-Infinity");
    boolean mns = f < 0; if (mns) f = -f;
    long lbits = Double.doubleToLongBits(f);
    int ue2 = (int)(lbits >>> 52 & 0x7ff);
    int e2 = ue2 - 1023 - 52 + (ue2==0 ? 1 : 0);
    long lquo, lmant = (lbits & ((1L << 52) - 1)) + (ue2==0 ? 0L : 1L << 52);
    int point = (int)Math.ceil(e2*llog2);
    BigInteger mant = BigInteger.valueOf(lmant);
    if (e2 > 0) {
        BigInteger num = mant.shiftLeft(e2 - point);
        lquo = roundQuotient(num, bipows5[point]);
        if (MantExpToDouble(lquo, point) != f) {
            num = num.shiftLeft(1);
            lquo = roundQuotient(num, bipows5[--point]);
            if (lmant==1L<<52 && MantExpToDouble(lquo, point) != f) {
                if (MantExpToDouble(++lquo, point) != f)
                    lquo = roundQuotient(num.shiftLeft(1), bipows5[--point]);
            }
        }
    } else {
        BigInteger num = mant.multiply(bipows5[-point]);
        BigInteger den = BigInteger.ONE.shiftLeft(point - e2);
        lquo = roundQuotient(num, den);
        if (MantExpToDouble(lquo, point) != f) {
            point--;
            num = num.multiply(BigInteger.TEN);
            lquo = roundQuotient(num, den);
            if (lmant==1L<<52 && MantExpToDouble(lquo, point) != f) {
                if (MantExpToDouble(++lquo, point) != f) {
                    point--;
                    lquo = roundQuotient(num.multiply(BigInteger.TEN), den);
                }
            }
        }
    }
    String sman = ""+lquo; int len = sman.length(), lent = len;
    while (sman.charAt(lent-1)=='0') {lent--;}
    StringBuilder str = new StringBuilder(23);
    if (mns) str.append('-');
    if (lent < 8 && point+len <= 0 && point+len > -3) {
        int zs = point+len; str.append("0."); while (zs++ < 0) str.append("0");
        return str.append(sman, 0, lent).toString();
    }
    if (lent < 8 && lent==point+len)
        return str.append(sman, 0, lent).append(".0").toString();
    str.append(sman, 0, 1).append('.').append(sman, 1, lent);
    if (lent==1) str.append('0');
    return str.append('E').append(point + len - 1).toString();
}
|endverbatim

%% \vfill
\beginsection{Performance}

\vbox{\settabs 2\columns
\+\hfill\figscale{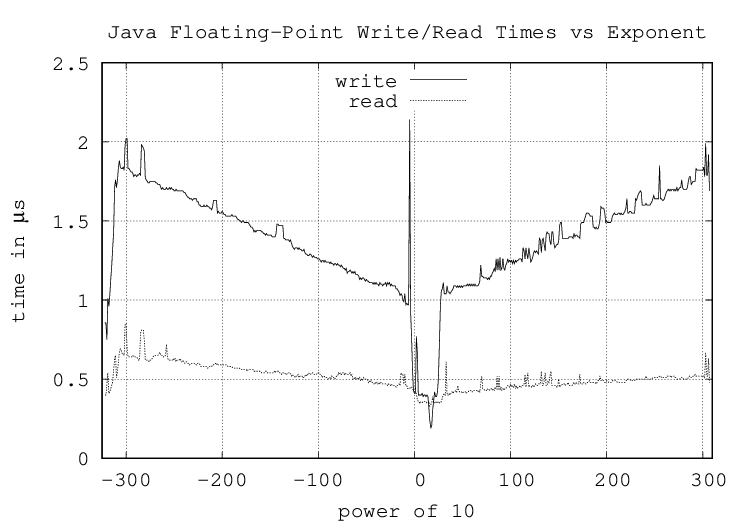}{225pt}\hfill&\hfill\figscale{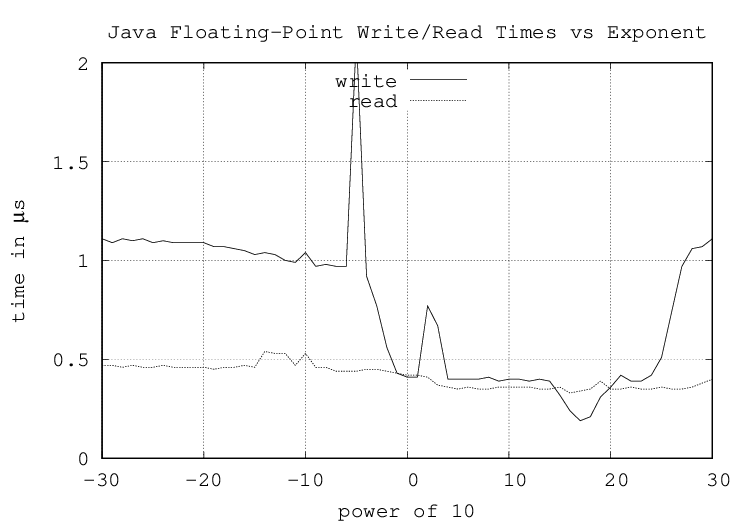}{225pt}\hfill&\cr
\+\hfill \figdef{rwtimesjna}\hfill&\hfill \figdef{rwtimesjnadet}\hfill&\cr
\+\hfill\figscale{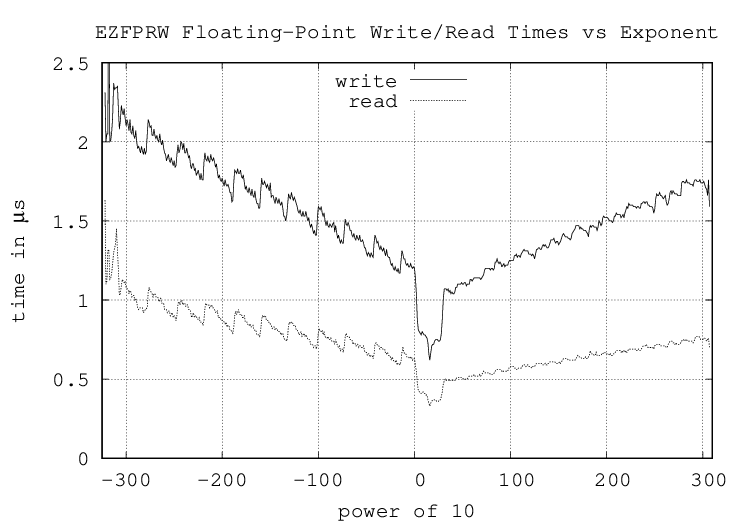}{225pt}\hfill&\hfill\figscale{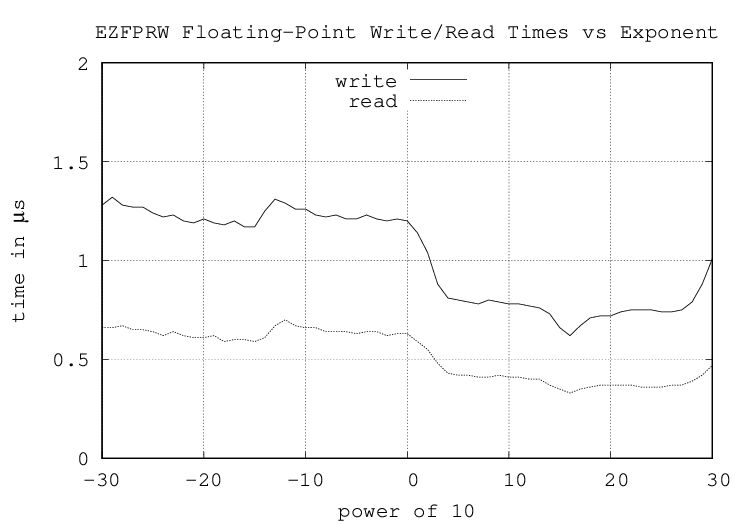}{225pt}\hfill&\cr
\+\hfill \figdef{rwtimesjbi}\hfill&\hfill \figdef{rwtimesjbidet}\hfill&\cr
\+\hfill\figscale{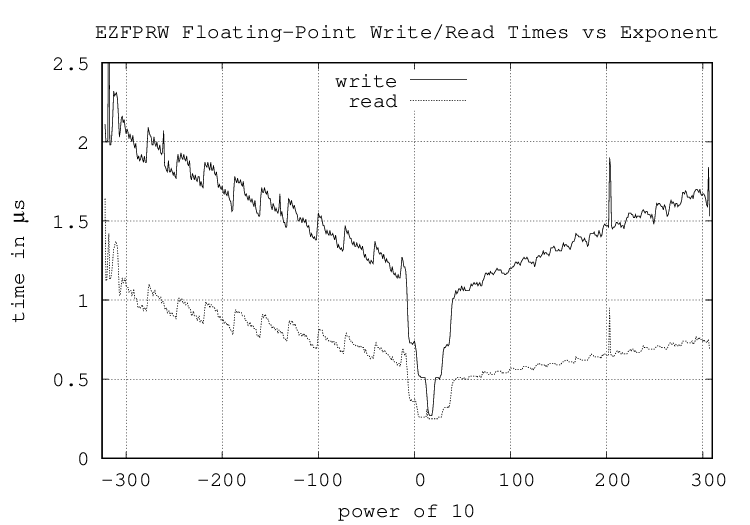}{225pt}\hfill&\hfill\figscale{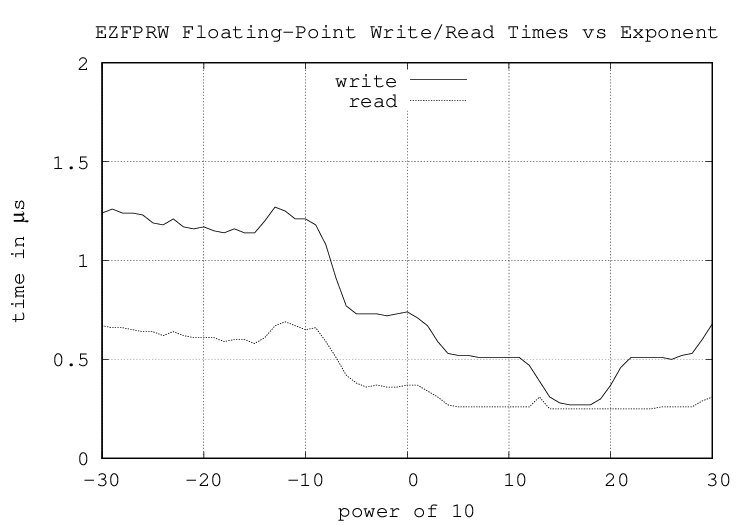}{225pt}\hfill&\cr
\+\hfill \figdef{rwtimesjbid}\hfill&\hfill \figdef{rwtimesjbiddet}\hfill&\cr
}

  IEEE-754 floating-point numbers have a finite range.  And the bulk
  of floating-point usage tends to have magnitudes within the range
  $1\times10^{-30}$ to $1\times10^{30}$.  Thus the asymptotic running
  time of floating-point conversion operations is of limited practical
  interest.  Instead, this article looks at measured running times of
  Java native conversions and conversions by these new algorithms over
  the full floating-point range.  These measurements were performed on
  Openjdk version "1.8.0\_171" running on a 2.40GHz Intel Core
  i7-5500U CPU with 16~GB of RAM hosting Ubuntu 16.04 GNU/Linux kernel
  4.4.0-127.

  A program was written which generated a vector of 100,000 numbers,
  $10^X$ where $X$ is a normally distributed random variable.  Then
  for each integer $-322\le n\le 307$, the vector of numbers is scaled
  by $10^n$, written to a file, read back in, and checked against the
  scaled vector.  The CPU time for writing and the time for reading
  were measured and plotted in \figref{rwtimesjna}.  An expanded view
  of \figref{rwtimesjna} for $-30\le n\le 30$ is plotted
  in \figref{rwtimesjnadet}.

%% \eject

  \figrefs{rwtimesjna} and \figrefn{rwtimesjnadet} show the performance of
  native conversions in Java version 1.8.0\_171.

  \figrefs{rwtimesjbi} and \figrefn{rwtimesjbidet} show the results
  for the power-of-5 {\tt BigInteger} algorithms implemented in Java.

  \figrefs{rwtimesjbid} and \figrefn{rwtimesjbiddet} show the results
  for the power-of-5 {\tt BigInteger} algorithms enhanced to use {\tt
  double}s and {\tt long}s instead of {\tt BigInteger} when the
  precision allows.

  The enhanced power-of-5 read algorithm is faster than Java native
  conversion in the exponent range $-5<n<30$.  Both write algorithms
  are at parity with Java native conversion for positive exponents;
  the read algorithm is roughly 50\% slower.  Denormalized
  conversions ($n<-309$) are less than half of the speed of Java
  native operations.  Over the rest of the negative exponent range the
  algorithms are roughly 50\% slower than Java native conversions.

  Because {\tt doubleToString} calls {\tt MantExpToDouble},
  improvements to the speed of {\tt MantExpToDouble} will benefit
  both.

  %% When the first (starved precision) conversion is skipped, {\tt
  %% doubleToString} is about as fast as Java native conversion in that
  %% exponent range as well; but, like Java-1.8 native conversions, some
  %% of the scientific-notation mantissas have one digit more than is
  %% necessary in order to read back the original floating-point value.

  %% Better performance than \figref{rwtimesjbidet} in the exponent range
  %% $10<n<20$ (with minimal precision mantissa) is achieved by doing the
  %% calculations in Java {\tt long} integers when they fit; the times
  %% are shown in \figrefs{rwtimesjbid} and \figrefn{rwtimesjbiddet}.

  %% There are four regions of the write and read curves
  %% of \figrefs{rwtimesjbi} and \figrefn{rwtimesjbid}.  In the range
  %% $0\le n\le 30$, the intermediate integers are small, fitting in a
  %% few cache lines.  For $n<-300$ the mantissa is unnormalized and
  %% requires smaller BigIntegers than for $n=-300$.  In the remaining
  %% regions the running time grows with the length of the intermediate
  %% BigIntegers, although the BigInteger computations take less than
  %% half of the overall conversion times.

\beginsection{Acknowledgments}

  Thanks to Tim Peters for finding a class of corner-cases which
  failed the first version of the {\tt doubleToString}.

\beginsection{Conclusion}

  The introduction of an integer {\tt roundQuotient} procedure
  facilitates algorithms for lossless (and minimal) conversions
  between (decimal) scientific-notation and (binary) IEEE-754
  double-precision floating-point numbers which are much simpler than
  algorithms previously published.

  Measurements of conversion times were conducted.  Implemented in
  Java, the optimized conversion algorithms executed faster than
  Java's native conversions over the range $10^{-5}$ to $10^{30}$, was
  comparable for writes of numbers greater than $10^{30}$, and 50\%
  slower over the rest of the IEEE-754 double-precision range.  The
  {\tt doubleToString} procedure is superior to Java native conversion
  in that it produces the minimum length mantissa which converts back
  to the original number.

\beginsection{References}

\bibliography{EZFPRW}
\bibliographystyle{alpha}

\vfill\eject
\bye